 \documentclass[a4paper,11pt]{article}
\usepackage{geometry,amssymb,amsmath,latexsym,epic,theorem}
\usepackage{xypic,xspace,enumerate}
\xyoption{curve}
\newcommand{\labto}[1]{\stackrel{#1}{\longrightarrow}}
{\theorembodyfont{\rmfamily}\newtheorem{example}{Example}[section]}
{\theorembodyfont{\rmfamily}}
\newtheorem{prop}[example]{Proposition}
\newtheorem{thm}[example]{Theorem}
{\theorembodyfont{\rmfamily}}
\newtheorem{cor}[example]{Corollary}

\newenvironment{proof}{\noindent {\bf Proof}}{\hfill $\Box$}

 \def\pt{\partial}
 \def\prt{\partial}

\def\eps{\varepsilon}
\def\leq{\leqslant}

\def\im{\mathrm{Im}\,}
\def\A{\alpha}

\def\cF{\mathcal{F}}
\def\cD{\mathcal{D}}
\def\Eq{\mathsf{Eq}}
\def\lan{\langle}
\def\ran{\rangle}

\newcommand{\threeaxes}[3]{\def\objectstyle{\scriptstyle}  \objectmargin={0pt}
\xy
(0,0)*+{}="a",(0,-6)*+{\rule{0em}{1.5ex}#2}="b",(7,0)*+{\;#1}="c",
(14,-3)*+{\;#3}="d" \ar@{->} "a";"b" \ar @{->}"a";"c"  \ar
@{->}"a";"d"\endxy }

\newcommand{\directs}[2]{\def\objectstyle{\scriptstyle}  \objectmargin={0pt}
\xy
(0,4)*+{}="a",(0,-2)*+{\rule{0em}{1.5ex}#2}="b",(7,4)*+{\;#1}="c"
\ar@{->} "a";"b" \ar @{->}"a";"c" \endxy }

\newcommand{\xdirects}[2]{\def\objectstyle{\scriptstyle}  \objectmargin={0pt}
\xy
(0,0)*+{}="a",(0,-6)*+{\rule{0em}{1.5ex}#2}="b",(7,0)*+{\;#1}="c"
\ar@{->} "a";"b" \ar @{->}"a";"c" \endxy }
\newcommand{\sdirects}[2]{\def\objectstyle{\scriptstyle}  \objectmargin={0pt}
\xy
(0,2.2)*+{}="a",(0,-2.5)*+{\rule{0em}{1.5ex}#2}="b",(7,2.2)*+{\;#1}="c"
\ar@{->} "a";"b" \ar @{->}"a";"c" \endxy }

\newcommand{\bl}{\mbox{\rule{0.08em}{1.7ex}\hspace{-0.00em}\rule{0.7em}{0.2ex}}}

\newcommand{\br}{\mbox{\rule{0.7em}{0.2ex}\hspace{-0.04em}\rule{0.08em}{1.7ex}}}

\newcommand{\tr}{\mbox{\rule[1.5ex]{0.7em}{0.2ex}\hspace{-0.03em}\rule{0.08em}{1.7ex}}}

\newcommand{\tl}{\mbox{\rule{0.08em}{1.7ex}\rule[1.54ex]{0.7em}{0.2ex}}}

\newcommand{\hh}{\mbox{\rule{0.7em}{0.2ex}\hspace{-0.7em}\rule[1.5ex]{0.70em}{0.2ex}}}

\newcommand{\vv}{\mbox{\rule{0.08em}{1.7ex}\hspace{0.6em}\rule{0.08em}{1.7ex}}}

\newcommand{\sq}{\mbox{\rule{0.08em}{1.7ex}\hspace{-0.00em}\rule{0.7em}{0.2ex}\hspace{-0.7em}\rule[1.54ex]{0.7em}{0.2ex}\hspace{-0.03em}\rule{0.08em}{1.7ex}}}

\def\epsilon{\varepsilon}
\def\C{\mathsf{C}}
\def\X{\mathsf{X}}
\def\Sets{\mathsf{Sets}}
\def\op{^\mathrm{op}}
\def\Grpd{\mathsf{Grpd}}
\def\Top{\mathsf{Top}}

\begin{document}
\title{Galois theory and a new homotopy double groupoid \\ of a map of spaces}

\author{Ronald Brown\thanks{Mathematics Division, School of Informatics, University of
Wales, Dean St., Bangor,  Gwynedd LL57 1UT,  U.K.  email:
r.brown{@}bangor.ac.uk}\\ \and George Janelidze\thanks{Mathematics
Institute, Georgian Academy of Sciences, Tbilisi, Georgia.  } }

\maketitle

\begin{center}
{\bf UWB Maths Preprint 02.18}
\end{center}

\begin{abstract}
The authors have used generalised Galois Theory  to construct a
homotopy double groupoid of a surjective fibration of Kan
simplicial sets. Here we apply this to construct a new homotopy
double groupoid of a map of spaces, which includes constructions
by others of a 2-groupoid,  cat$^1$-group or crossed module. An
advantage of our construction is that the double groupoid can give
an algebraic model of a foliated bundle.\footnote{2000 Maths
Subject Classification: 18D05, 20L05, 55 Q05, 55Q35}
\end{abstract}

\section*{Introduction}

Our aim is to develop for  any map $q:M \to B$ of topological
spaces the construction and properties of a new homotopy double
groupoid which has the form of the left hand square in the
following diagram, while the right hand square gives a morphism of
groupoids:
\begin{equation}\label{Rho}\xymatrix @=3pc {\rho_2(q) \ar @<1ex> [r] ^s \ar @<-1ex> [r]
_t \ar @<1ex> [d]  \ar @<-1ex> [d] & \pi_1(M)   \ar[l]  \ar @<1ex>
[d]
 \ar @<-1ex> [d]\ar [r] ^{\pi_1q}& \pi_1(B) \ar @<1ex> [d]  \ar @<-1ex> [d] \\
\mathsf{Eq}(q) \ar [u]  \ar @<1ex> [r] ^s \ar @<-1ex> [r] _t & M
\ar [l] \ar[u]  \ar [r] _q & B\ar [u]}\end{equation} where:\\
$\pi_1(M)$ is the fundamental groupoid of $M$;\\  $\mathsf{Eq}(q)$
is the equivalence relation determined by $q$; and\\  $s,t$ are
the source and target maps of the groupoids. \\ Note that $qs=qt$
and $(\pi_1q)s=(\pi_1q)t$, so that $\rho_2(q)$ is seen as a double
groupoid analogue of $\Eq(q)$.

 This double groupoid contains the
2-groupoid associated to a map defined by Kamps and Porter in
\cite{KP}, and hence also includes the cat$^1$-group of a
fibration defined by Loday in \cite{Lo}, the 2-groupoid of a pair
defined by Moerdijk and Svensson in \cite{MandS:2-types}, and the
classical fundamental crossed module of a pair of pointed spaces
defined by J.H.C. Whitehead. Advantages of our construction are:
\\ (i)   it contains information on the map $q$, and \\ (ii) we get
different results if the topology of $M$ is varied to a finer
topology. \\ In particular, we can apply this construction in the
case $M$ is foliated by replacing the topology on $ M$ by a finer
one so that $\pi_1M$ is replaced by the fundamental groupoid of
the foliation.

This applies in particular to the M\"{o}bius Band with its
standard foliation by circles. We can extract from this double
groupoid a small version  $\cD(M)$ with only three vertices, and
which seems to represent well many properties of the M\"{o}bius
Band. It has basic vertices, edges and squares as follows:
$$\def\labelstyle{\textstyle}
\xymatrix{ B  \ar [r] ^\theta \ar @{} [dr]|\alpha \ar [d] _\eta &
C \ar [d]^\xi \\
A \ar [r] |\iota & A \\
C \ar [u] ^\xi \ar [r] _\phi \ar @{} [ur] |\beta & B \ar [u]_\eta
} \qquad \xdirects{2}{1}
$$
Note that the vertical groupoid for $\circ_1$ is an indiscrete
groupoid, while the horizontal groupoid for $\circ_2$ contains a
copy of the infinite cyclic group, since there are compositions
$$a_1\circ_2a_2\circ_2 \cdots \circ_2a_n$$ where the $a_i$ are
alternately $(-_1\alpha)$ and $ \beta$.

The idea for this double groupoid arose from the Generalised
Galois Theory of Janelidze \cite{J,J2}, which under certain
conditions gives a Galois groupoid from a pair of adjoint
functors. The standard fundamental group arises from the adjoint
pair between topological spaces and sets given by discrete and
$\pi_0$, see for example \cite{BJ1}. The adjoint pair between
simplicial sets and crossed complexes given by nerve and $\pi_1$
was studied in \cite{BJ2} and shown to lead to a Galois double
groupoid of a fibration of simplicial sets. We are now giving a
topological version of this construction. We show that if $p:E \to
B$ is a Serre fibration then the fundamental groupoid $\pi_1(E)$
has an additional compatible groupoid structure arising from the
equivalence relation $\Eq(p)$ defined by the map $p$; these two
groupoid structures define a double groupoid which we write
$\gamma(p)$, since it is defined by methods of Galois theory. The
double groupoid $\rho(q)$ arises by pullback by $i$ where $q=pi$
is the usual factorisation of any map through a homotopy
equivalence $i$ and a fibration $p$. However the proof of the
relation of $\rho(q)$ with classical notions and compositions is
tricky, and so is given in some detail. A further reason for this
detail is the possibility that a modification of this construction
could be used in association with the `thin fundamental groupoids'
and their smooth structures in differential geometrical
situations, as exemplified by Mackaay and Picken in \cite{MP}.

Here is some background to the search for higher groupoid models
of homotopical structures (for more detailed references, see
\cite{Br}). Geometers in the early part of the 20th century were
aware that in the connected case the first homology group was the
fundamental group made abelian, and that homology groups existed
in all positive dimensions. Further, the fundamental group gave
more information in geometric and analytic contexts than did the
first homology group. They were therefore interested in seeking
higher dimensional versions of the non abelian fundamental group.
E. \v{C}ech submitted to the 1932 ICM at Zurich a paper on higher
homotopy groups, using maps of spheres. However these groups were
quickly proved to be abelian in dimensions $> 1$, and on this
ground \v{C}ech was asked to withdraw his paper, so that only a
small paragraph appeared \cite{Ce}. Thus the dream of these
topologists seemed to fail, and was widely felt to be a mirage,
although the abelian higher homotopy groups became and still are
very important.

J.H.C. Whitehead in the 1940s introduced the notion of crossed
module, using the boundary of the second relative homotopy group
of a pair and the action of the fundamental group. He and Mac~Lane
showed that crossed modules classified (connected) homotopy
2-types. Crossed modules are indeed more complicated than groups,
and they make a good candidate for `2-dimensional groups'.

In the 1960s, Brown introduced the fundamental groupoid of a space
on a set of base points, and the writing of his 1968 book on
topology suggested to him that all of 1-dimensional homotopy
theory was better expressed in terms of groupoids rather than
groups. This raised the question of the putative value of
groupoids in higher homotopy theory. A relation of certain double
groupoids to crossed modules was worked out with C.B. Spencer in
the early 1970s, and this showed that double groupoids are indeed
more complicated than groups. A definition of a homotopy double
groupoid of a pair of pointed spaces was made with P.J. Higgins in
1974, and exploited to obtain a 2-dimensional Van Kampen type
theorem for this double groupoid, and hence for Whitehead's
crossed module of a pair (see \cite{BH1}). The double groupoid
constructed in \cite{BH1} is edge symmetric and has a connection,
and so is not the same as that constructed  here.

A classification of certain double groupoids is given in
\cite{BM}, but this does not yield much information for the double
groupoid considered here. Thus there is still a way to go in the
understanding and in the use of double groupoids.

Higher homotopy groupoids were defined by Brown and Higgins for a
filtered space in \cite{BH3}, and by Loday for an $n$-cube of
spaces in \cite{Lo}; his cat$^n$-groups were shown there to model
connected homotopy $(n+1)$-types. These higher groupoid methods
yield new calculations in homotopy theory through higher order Van
Kampen theorems \cite{BH3,BL}, as well as suggesting new algebraic
constructions.

\section{Galois groupoids}

 Later we will be considering the category $\C =
\Sets^{\Delta^{op}}$ of simplicial sets and the fundamental
groupoid functor $I=\pi_1 : \C\to \X$  from the category $\C$ to
the category $\X = \Grpd$ of (small) groupoids. Further,  $C$, an
internal category in $\C$,  will be  the particular simplicial
category (actually groupoid) $\Eq(p)$ which is the equivalence
relation (in $\C$) determined by $p$ where $p : E\to B$ is a
surjective fibration of Kan complexes. Here we give first the
general result, using this notation.

Let $I : \C \to  \X$ be an arbitrary functor between categories
$\C$ and $\X$ with pullbacks, and let
\begin{equation}C=\left(\xymatrix{C_2 \ar [r] \ar@<1ex>[r] \ar@<-1ex>[r]
& C_1\ar@<1ex>[r] \ar@<-1ex>[r]& C_0 \ar [l] }\right)
\end{equation} be an internal category in $\C$. We recall

\begin{prop} Suppose the canonical morphisms
\begin{align} &I(C_1\times_{C_0} C_1)\to
 I(C_1)\times _{I(C_0)}I(C_1)\\
          & I(C_1 \times_{C_0} C_1\times_{C_0}C_1)\to
           I(C_1)\times_{I(C_0)}I(C_1)\times_{I(C_0)}I(C_1)
          \end{align}
are isomorphisms. Then:\begin{enumerate}[(a)] \item
\begin{equation} I(C) =\left(\xymatrix{I(C_2)\ar@<1ex>[r]
 \ar@<-1ex>[r]\ar [r]&
I(C_1)\ar@<1ex>[r] \ar@<-1ex>[r]& I(C_0) \ar [l] }\right)
\end{equation} is an internal category in $\X$; \item if $C$ is a
groupoid, then so is $I(C)$. \end{enumerate}
\end{prop}

For a morphism $p : E\to  B$ in $\C$, let $\Eq(p) =$
\begin{equation}\xymatrix{
(E\times_{B}E)\times_{E}(E\times_{B}E)\approx
E\times_{B}E\times_{B}E\ar [r] \ar@<1ex>[r]  \ar@<-1ex>[r] &
E\times_{B}E  \ar@<1ex>[r] \ar@<-1ex>[r] & \ar [l]      E}
\end{equation}
be the equivalence relation corresponding to $p$ (=kernel pair of
$p$) regarded as an internal groupoid in $\C$. Applying
Proposition 1.1 with $C = \Eq(p)$, we obtain:
\begin{cor}Suppose the canonical morphisms
\begin{align} &I((E\times_{B} E)\times_{E}(E\times_{B}E))\to
 I(E\times_{B} E)\times _{I(E)}I(E\times_{B} E)\\
          &I(E\times_{B} E) \times_{E}(E\times_{B} E)\times_{E}(E\times_{B} E))\to
           I(E\times_{B} E)\times_{I(E)}I(E\times_{B} E)\times_{I(E)}I(E\times_{B} E)
           \end{align}
are isomorphisms. Then $I(\Eq(p)) =$
\begin{equation}\xymatrix{
I((E\times_{B}E)\times_{E}(E\times_{B}E))\approx
I(E\times_{B}E\times_{B}E)\ar [r] \ar@<1ex>[r]  \ar@<-1ex>[r] &
I(E\times_{B}E)  \ar@<1ex>[r] \ar@<-1ex>[r] & \ar [l] I(E)}
\end{equation}
is an internal groupoid in $\X$.
\end{cor}
This fact, which goes back to A. Grothendieck's observation ``the
fundamental groupoids are to be defined as quotients of
equivalence relations", is used in categorical Galois theory and
its various special cases (see [1], [4], and other references in
[1]), to define the {\em Galois groupoid} of $(E,p)$ as
\begin{equation}
          Gal_I(E,p) = I(\Eq(p)).
\end{equation}
In particular this applies to the following situation studied by
the authors before (see Proposition 3.5 in [2]):
\begin{prop}
 Let $I : \C\to  \X$ be the fundamental groupoid
functor from the category $\C = \Sets^{\Delta^{op}}$ of simplicial
sets to the category $\X = \Grpd$  of (small) groupoids, and $p :
E\to B$ a surjective fibration of Kan complexes. Then the
morphisms (6) and (7) are isomorphisms and so the Galois groupoid
(9) is well defined. Since the internal groupoids in $\Grpd$ are
the same as double groupoids, it is a double groupoid.\end{prop}

\section{From simplicial sets to topological spaces}

 Consider the diagram
 \begin{equation}
\xymatrix{ \Top \ar@<-0.5ex> [r]_-S &\ar @<-0.5ex>[l]_-R
\Sets^{\Delta\op} \ar
@<0.5ex>[r] ^-I & \Grpd\ar @<0.5ex>[l]^-H\\
& \Delta  \ar [ul]^r \ar [u] _y \ar [ur]_i}
 \end{equation}
 where
 \begin{itemize}
\item $\Top$ is the category of topological spaces, $R$ is the
geometric realisation functor, and $S$ is its right  adjoint,
usually called the singular complex functor; \item $ I \vdash  H$
is the adjoint pair used in [2], i.e. $I$ is the fundamental
groupoid functor, and $H$ the nerve functor; \item $y$ is the
Yoneda embedding, $ r$ and $i$ are  the restrictions of $R$ and
$I$ respectively along $y$; explicitly, $r$ is the singular
simplex functor and $i$ carries finite ordinals to codiscrete
groupoids on the same sets of objects.
\end{itemize}

By the universal property of the Yoneda embedding, the two
adjunctions of the row are uniquely (up to isomorphisms)
determined by $r$ and $i$; let us also recall from \cite{BJ2} and
\cite{J,J2}:

\begin{prop} (a) The composite $IS : \Top\to  \Grpd$ can be identified with
the classical (geometric) fundamental groupoid functor $\pi_1$;

(b) For every topological space $X$, $S(X)$ is a Kan complex;

(c) The $S$-image of a morphism $p$ in $\Top$ is a Kan fibration
if and only if $p$ itself is a Serre fibration.\end{prop}

\section{What is the Galois double groupoid of a Serre
fibration?}
 Let $p : E\to B$ be a Serre fibration of topological
spaces. By Propositions 1.3 and 2.1, the Galois (double) groupoid
$Gal_I(S(E),S(p))$ is well defined. Moreover, since the functor
$S$ being a right adjoint preserves pullbacks, we can write
\begin{equation}
          Gal_I(S(E),S(p)) \approx Gal_{SI}(E,p)\approx  Gal_{\pi_1}(E,p)
\end{equation}
and conclude that $Gal_{\pi_1}(E,p)$ also is a well-defined double
groupoid. We  write this double groupoid as $\gamma(p)$ to
indicate the relation with Galois theory, and now describe it
explicitly.

The underlying double graph has the following description:
\begin{itemize}
\item presented as an internal groupoid in $\Grpd$, $\gamma(p)$
is displayed as
\begin{equation}
\xymatrix{ \pi_1((E\times_{B}E)\times_{E}(E\times_{B}E))\approx
\pi_1(E\times_{B}E\times_{B}E)\ar [r] \ar@<1ex>[r]  \ar@<-1ex>[r]
& \pi_1(E\times_{B}E)  \ar@<1ex>[r] \ar@<-1ex>[r] & \ar [l]
\pi_1(E)}\end{equation}
 from which we conclude:
\item  the set of objects in $\gamma(p)=Gal_{\pi_1}(E,p)$ is $E$;
\item a horizontal arrow $e\to e'$ is a morphism $e\to  e'$ in
$\pi_1(E)$, i.e. a homotopy class of a path from $e$ to $e'$
(recall that such a homotopy $h : [0,1]\times[0,1]\to E$ is
required to have $h(0,t) = e$ and $h(1,t) = e'$ for every $t$ in
$[0,1]$); \item  a vertical  arrow $e\to e'$ is just the pair
$(e,e')$ provided $p(e) = p(e')$; \item a square
\begin{equation}\def\labelstyle{\textstyle}
\xymatrix{ e_1\ar[r] ^\phi \ar [d] \ar @{}[dr] |u & e_2 \ar [d]  \\
e_1' \ar [r] _{\phi'} & e_2'}
\end{equation}
 is a homotopy class of a path from $(e_1,e'_1)$ to $(e_2,e'_2)$
 in $E\times_BE$; its vertical  domain $\phi$  and vertical
 codomain $\phi'$ are homotopy classes of paths from $e_1$ to $e_2$ and
 from $e'_1$ to $e'_2$ respectively, and the horizontal are pairs
 $(e_1,e'_1)$
 and $(e_2,e'_2)$ respectively.

 \end{itemize}

Clearly there is no problem with the horizontal composition as
well - since we know that $\pi_1(E\times_BE)$ and $\pi_1(E)$ are
groupoids. The only non-trivial part of our construction is the
vertical composition of squares:

Given
\begin{equation}\def\labelstyle{\textstyle}
\xymatrix{e_1 \ar [r] ^\phi \ar [d] \ar @{}[dr] |u & e_2 \ar [d]  \\
e_1' \ar [r] |{\phi'}\ar @{}[dr] |{u'}\ar[d]  & e_2'\ar [d]\\
e_1'' \ar [r] _{\phi'} & e_2'' } \quad \xdirects{2}{1}
\end{equation}
we have to define $u \circ_1 u'$, which must be (an equivalence
class of) a path from $(e_1,e''_1)$ to $(e_2,e''_2)$ in
$E\times_BE$. Of course we should choose a representative $f$ of
$u$, a representative $f'$ of $u'$, a homotopy $h$ between the
vertical codomain path of $f$ and the vertical domain path of
$f'$, paste them together, and take the homotopy class of the
resulting path in $E\times_BE$. However, how do we know that such
an $u\circ_1 u'$ does not depend on the choices we made? The nice
consequence of the results above is that we do not need to prove
this. Indeed, since the morphism
$\pi_1((E\times_BE)\times_E(E\times_BE))\to
\pi_1(E\times_BE)\times_{\pi_1(E)}\pi_1(E\times_BE)$ is an
isomorphism, the pair $(u',u)$ determines a path $v$ in
$(E\times_BE)\times_E(E\times_BE)$, and the desired vertical
composite $u'u$ is nothing but the image of $v$ under the functor
$\pi_1((E\times_BE)\times_E(E\times_BE))\to \pi_1(E\times_BE)$
induced by the composition map
$(E\times_BE)\times_E(E\times_BE)\to E\times_BE$ of $\Eq(p)$.
Moreover, we do not have to worry about associativity of the
horizontal composition, which in fact follows from
$$\pi_1((E\times_BE)\times_E(E\times_BE)\times_E(E\times_BE)) \to
\pi_1(E\times_BE)\times_{\pi_1(E)}\pi_1(E\times_BE)\times_{\pi_1(E)}\pi_1(E\times_BE)$$
being an isomorphism.

\begin{example} \label{prodfib}
Suppose for example that $p: E=( B \times F) \to B$ is the
projection of a product, and so is a fibration. Then $\Eq(p)= E
\times _B E$ is homeomorphic to the product  $B \times (F \times
F)$ and hence $\pi_1(\Eq(p))$ is easily determined, together with
its double groupoid structure. The vertical edges are triples
$(b,f^-,f^+), b \in B, f^\pm \in F$, the  horizontal edges are
pairs $(\beta,\phi) \in \pi_1 B \times \pi_1 F $ and the squares
are triples $$(\beta, \phi^-,\phi^+) \in \pi_1(B) \times \pi_1(F)
\times \pi_1(F)$$ with vertical boundaries given by $\pt^\pm_1
(\beta, \phi^-,\phi^+)=(\beta,\phi^\pm)$. The horizontal
composition $\circ_2$ is that of the fundamental groupoids; the
vertical composition $\circ_1$ reflects that of the product
groupoid $B \times (F \times F)$ in which $B$ is the discrete
groupoid and $F \times F$ is the codiscrete groupoid.

\end{example}

In the  next sections we translate the construction for
$\gamma(p)$ where $p$ is a fibration into results for an arbitrary
map $q: M \to B$ by the usual factorisation process, and so relate
these ideas to classical constructions.
\section{The homotopy double groupoid of an arbitrary continuous map}

If $C$ is an internal category (or a groupoid) in a category $\X$
with pullbacks, and $i : M\to  C_0$ a morphism into the
object-of-object of $C$, we can always form the ``induced internal
category (respectively groupoid)" $i^*(C) =$
\begin{equation}
\xymatrix{C_2\times_{C_0 \times C_0 \times C_0} (M \times M \times
M) \ar [r] \ar@<1ex>[r] \ar@<-1ex>[r] & C_1\times_{C_0 \times
C_0}(M \times M)\ar@<1ex>[r] \ar@<-1ex>[r]& M \ar [l] }
\end{equation}
which, in the case $\X = \Sets$, can be simply described as the
category with objects all elements of $M$, and morphisms ``as in
$C$". In particular, for a topological space $E$ and a subspace
$M$ one defines the relative fundamental groupoid $\pi_1(E,M)$ as
$i^*(\pi_1(E))$, where $i : M\to  E$ is the inclusion map; in this
paper we do not consider $\pi_1(E)$ as a topological groupoid, but
the same construction of $\pi_1(E,M)$ can be repeated internally
to $\Top$.

Let $p : E\to  B$ be a Serre fibration and $i : M\to  E$ an
arbitrary map in $\Top$. Using the morphism $\pi_1(i) :
\pi_1(M)\to \pi_1(E)$ of fundamental groupoids, and having in mind
that $\pi_1(E)$ is the object-of-objects of $Gal_{\pi_1}(E,p)$
(when it is regarded as an internal groupoid in $\Grpd$), we can
construct what we are going to call the relative Galois double
groupoid $Gal_{\pi_1}(E,p,M,i)$ as
\begin{equation}
          Gal_{\pi_1}(E,p,M,i) = \pi_1(i)^*Gal_{\pi_1}(E,p).    \end{equation}
A good reason for introducing this new double groupoid is that any
continuous map $q: M\to  B$ can be presented as a composite $q =
pi$ above with $i$ being a homotopy equivalence, and so any $q$
determines such a double groupoid equivalent (in an appropriate
sense) to the Galois double groupoid of a Serre fibration. It is
this double groupoid we have displayed in diagram (1).

In the next sections we will identify this double groupoid in
terms of  homotopy classes of certain maps and so relate these
ideas and the compositions obtained more clearly in terms of
classical homotopical ideas. One further reason for doing this is
to allow the possibility of generalisations of these geometric
constructions to situations where the appropriate Galois theory is
not yet obtained, for example in terms of smooth structures and
`thin fundamental groupoids', as in \cite{MP}. This could lead to
smooth structures on variations of the constructions given here.

\section{Geometric interpretation}
We now start to interpret the previous results geometrically in
terms of compositions related to classical notions of relative
homotopy groups. To this end, we use standard notation for the
cubical singular complex $KM$ of a space $M$. Here $K_nM$ is the
set of singular $n$-cubes in $M$ (i.e. continuous maps $I^n \to M$
with $K_0M$ identified with $M$). There are standard {\em face
maps}
 $\pt^\pm_i: K_nM
\to K_{n-1}M$ and {\em degeneracy maps} $\eps_i :K_{n-1}M \to
K_{n}M$ for $i=1, \ldots,n$ (with $\eps_1:K_0M \to K_1M$ written
simply $\eps$, and giving for $x \in M$ the constant path $\eps x$
at $x$). There are also for $i=1,\ldots,n$ {\em compositions}
$a\circ_ib$ defined for $a,b\in K_nM$ such that
$\pt^+_ia=\pt^-_ib$, and {\em inversions} $-_i: K_n \to K_n$.
Finally we shall later use {\em connections} $ \Gamma_i^\A :K_{n-1
}M \to K_nM$ for $i=1,\ldots,n,\A=\pm$ induced by the maps
$\gamma_i^\A : I^n \to
 I^{n-1}$    defined by
   $$   \gamma _i^\A (t_1 ,t_2 ,\ldots,t_n ) =
             (t_1 ,t_2 ,\ldots,t_{i-1},A(t_i ,t_{i+1}),t_{i+2},\ldots,t_n )
             $$
 where $A(s,t)=\max(s,t), \min(s,t)$ as $\A=-,+$ respectively.
The enrichment with connections $\Gamma^+_i$ for the traditional
cubical sets was introduced in \cite{BH2}.  The full properties of
these structures are set out in for example \cite{AABS}. Here we
will assume only the obvious geometric properties in the range
$n=0,\ldots,3$.

Let $q: M \to B$ be a map of spaces. We recall the standard
factorisation $q=pi$ where $i: M \to M_q$ is a homotopy
equivalence and $p:M_q \to B$ is a fibration. Here
$$M_q=\{(x,\lambda)\mid \lambda:I \to B, \,\lambda(0)=q(x)\}
\subseteq M \times B ^I$$ and $p(x,\lambda)=\lambda(1)$, while $i:
x \mapsto (x,\eps(q(x))$ where $\eps(q(x))$ is the constant path
at $q(x)$ in $B$. A point of $M_q \times _B M_q$ is a  pair
$((x,\lambda),(x',\lambda'))$ with $\lambda(1)=\lambda'(1)$ as in
the following picture:
$$\def\labelstyle{\textstyle}\xymatrixcolsep{0.3pc}
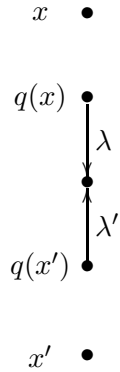
\xymatrix@M=0pt{x & \bullet \\
q(x)& \bullet \ar [d] ^\lambda \\
&\bullet\\q(x')& \bullet \ar [u]_{\lambda'} \\x'&\bullet}$$ So
$M_q \times _BM_q$ is homeomorphic to the space $M'$ of triples
$(x,\mu,x') \in M \times K_1B \times M$ such that $\mu(0)=q(x),
\mu(1)=q(x')$. Hence we have a groupoid structure on $M'$ with
object set $M_q$, where the source and target maps $s',t'$ send
$(x,\mu,x')$ to $(x,\mu_1)$, $(x',-\mu_2)$ where $\mu_1,\mu_2 \in
K_1B$ are respectively the rescaled forms of the first half and
the second half of $\mu$. The composition $\circ'$ in this
groupoid is, when defined, given by
$$(x,\mu,x')\circ'(x',\mu',x'')=(x',\mu_1 \circ_1 \mu_2,x'')$$
where $\circ_1$ is here the usual composition of paths. Of course
this composition is defined if and only if $\mu_2=-_1\mu_1$.

Let $j: M \times_BM \to M'$ be given by $(x,x') \mapsto (x,
\eps(qx),x')$. The set of paths $I \to M'$ with end points in $\im
\, j$ can be identified with the subset $R_2(q)$ of  $K_1M \times
K_2B \times K_1M$ of triples $(f,\A,g)$ such that $qf=\prt^-_1 \A
, qg=\prt^+_1 \A$ and $\prt^-_2 \A,\prt^+_2 \A$ are constant
paths. Thus an element of $R_2(q)$ may be pictured as:
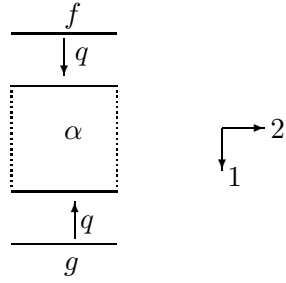
\begin{center}
\setlength{\unitlength}{0.7mm}
\begin{picture}(40,40)(0,-8)

\put(0,-10){\line(1,0){20}} \put(0,30){\line(1,0){20}}
\put(10,32){$f$}
 \put(10,29){\vector(0,-1){7}}
 \put(12,25){$q$}
 \put(12,-9){\vector(0,1){7}}
 \put(13,-7){$q$}
 \put(10,-15){$g$}
  \put(10,10){$\alpha$}
  \put(0,0){\line(1,0){20}}
  \put(0,20){\line(1,0){20}}
 \put(40,12){\vector(1,0){8}} \put(40,12){\vector(0,-1){8}}
 \put(49.2,10.2){2} \put(41,1){1}
  {\dottedline{1}(0,1)(0,20)}
{\dottedline{1}(20,0)(20,20)}

\end{picture}
\end{center}
\vspace{1ex} where the dotted lines show constant paths. Thus
$R_2(q)$ fits in the following diagram:
\begin{equation}\label{R}\xymatrix @=3pc{R_2(q) \ar @<1ex> [r]  \ar @<-1ex> [r]\ar @<1ex>
[d]  \ar @<1ex> [d]\ar @<-1ex> [d] & \ar @<1ex> [d]\ar @<-1ex> [d] \;\; K_1M   \ar[l]\\
\mathsf{Eq}(q) \ar [u]  \ar @<1ex> [r]  \ar @<-1ex> [r]  & M \ar
[l] \ar[u] } \quad \xymatrix@=1.5pc{&\\& \xdirects{2}{1} \\&
\;\\\;&}
\end{equation}
The boundary maps are given by:
\begin{align*}
\prt^-_1(f,\A,g) &= f, \\ \prt^+_1(f,\A,g) &= g, \\
\prt^-_2(f,\A,g) &= (f(0), g(0)), \\ \prt^+_2(f,\A,g) &= (f(1),
g(1)). \\ \intertext{The degeneracy maps $\eps_1: K_1M \to R_2(q),
\; \eps_2: \mathsf{Eq}(q) \to R_2(q)$ are given by:} \eps_1(f)&= (f, \eps_1(qf),f), \\
\eps_2(x,x')&= (\eps(x), \eps_1^2(q(x)), \eps(x')),
\end{align*}
for $(f,\A,g) \in R_2(q), (x,x') \in M \times_B M.$

Clearly \eqref{R} can be considered as a diagram of reflexive
graphs. We now examine compositions on $R_2(q)$.

 The set $R_2(q)$
 has two partial compositions. The composition $\circ_2$ is
 determined
by the usual composition of paths and squares in this direction:
$$(f,\A,g) \circ_2 (f',\A',g')= (f \circ_1f', \A\circ_2\A',g\circ_1g').$$
The composition $\circ_1$ in the direction 1 is given by
\begin{align}
(f,\A,g) \circ_1 (g,\beta, h) &= (f, \A\circ_1\beta,h).
\end{align}

Note that  this definition generalises  a construction by Kamps
and Porter in \cite[Section 4.1]{KP},  in which they assume
$f(0)=g(0), f(1)=g(1)$ whereas in our situation we have only
$qf(0)=qg(0), qf(1)=qg(1)$. Hence they end up with a 2-groupoid,
and we end up with a double groupoid. Their method of proving the
properties of their homotopy 2-groupoid is to assume first that
$p$ is a fibration, and then apply this case to an arbitrary map
$q$ by converting it to a fibration $p= \bar{q}$. This is
analogous to our methods, except that we have used Galois theory,
whereas they use directly properties of fibrations.

We now form the quotient of diagram \eqref{R} by taking homotopy
classes rel vertices of $R_2(q)$ and of $K_1M$ to yield the
diagram:
\begin{equation}\label{Rho2}\xymatrix @=3pc {\rho_2(q) \ar @<1ex> [r]  \ar @<-1ex> [r]
 \ar @<1ex> [d]  \ar @<-1ex> [d] & \pi_1(M)   \ar[l]  \ar @<1ex> [d]
 \ar @<-1ex> [d]\\
\mathsf{Eq}(q) \ar [u]  \ar @<1ex> [r]  \ar @<-1ex> [r]  & M \ar
[l] \ar[u] }\end{equation} where $\pi_1(M)$ is the fundamental
groupoid of $M$. It is clear that the horizontal composition
$\circ_2$ on $R_2(q)$ is inherited by $\rho_2(q)$. Our main result
of this section is a direct verification that the composition
$\circ_1$ is also inherited, without going through the simplicial
Galois theory of the previous sections. We also have to show that
this composition is related to that derived from the equivalence
relation structure on $K_1M'$.

In fact we prove a stronger result. The  set $$ R_2(q) {\;}
_{\prt^+_1}\!\! \times _{\prt^-_1}R_2(p)$$  is the domain of
composition of $\circ_1$ on $R_2(q)$. A {\em homotopy rel
vertices} on this set is a continuous family
$((f_u,\A_u,g_u),(g_u,\beta_u, k_u)), 0 \leq u \leq 1$ of elements
of this set such that $f_u(0)=f_0(0), k_u(1)=k_0(1), \; 0\leq u
\leq 1$. We use the notation $\pi_0^v$ for the set of homotopy
classes rel vertices.
\begin{thm}\label{bijection}
The natural map
\begin{equation}
  \Xi:   \pi_0^v( R_2(q) {\,}_{\prt^+_1}\!\! \times _{\prt^-_1}R_2(q)) \to
            \rho_2(q) {\,}_{\prt^+_1}\!\! \times _{\prt^-_1}\rho_2(q)
\end{equation}
defined by the projections, is a bijection.
\end{thm}

For the proof we use properties of the connections, and we use the
following notation from \cite{Spe}.

We write:

$$ \tl \quad  \sdirects{j+1}{j} \text{ for } \Gamma^+_j;
\quad \br \quad  \sdirects{j+1}{j} \text{ for }\Gamma^-_j$$
$$ \vv \quad  \sdirects{j+1}{j} \text{ for } \eps_{j};
\quad \hh \quad  \sdirects{j+1}{j} \text{ for }\eps_{j+1}.$$ Thus
the thick lines denote degenerate faces. We shall use inversions
applied to connections, for example
$$-_2 \tl , \quad \sdirects{2}{1}, $$ and write this also as $\tr\;$
since it coincides with $-_1\br$.

This notation allows us to write some compositions  as for example
that involving 3-cubes $x,y$ with $\pt^+_3 x = \pt ^-_3 y$ as
$$ A = \begin{bmatrix} \vv & \tl \; \\ x & y  \end{bmatrix}\quad  \directs{3}{2} $$
which is an abbreviation for
$$  \begin{bmatrix} \eps_{2}\pt^-_2x  & \Gamma^+_2\pt^-_2y \; \\ x & y
\end{bmatrix}$$ and
makes it transparent what are the faces of $A$. The direction
arrows are omitted when convenient.

\begin{proof}{\bf \  of theorem \ref{bijection}}
We define an inverse $\Phi$ for $\Xi$.

We use square brackets $[\;]$ to denote homotopy classes. Let
$([f, \A,g], [h,\beta,k]) \in  \rho_2(q) {\,}_{\prt^+_1}\!\!
\times _{\prt^-_1}\rho_2(q) $. Then there is a homotopy rel
vertices of paths $\xi: g \simeq h: I^2 \to M$. We set
\begin{equation} \label{P}
\Phi ([f, \A,g], [h,\beta,k]) = [(f,\A,g),(g,(q \xi)\circ_1 \beta,
k)]
\end{equation}
and have to prove $\Phi$ is well defined and an inverse to $\Xi$.

Suppose we are given homotopies
\begin{align}
\kappa: (f, \A,g) &\equiv (f',\A',g') \\
\lambda: (h, \beta,k) & \equiv (h', \beta ',k')\\
\xi': g'&\simeq h'.
\end{align}
Then $\kappa, \lambda$ are given by three component homotopies rel
vertices
\begin{gather}
\kappa_1: f \simeq f', \; \kappa_2: \A \equiv \A', \kappa_3: g
\simeq g', \\ \lambda_1:h \simeq h', \; \lambda_2: \beta \equiv
\beta', \lambda_3: k \simeq k',
\end{gather} with the properties that
\begin{gather}q\kappa_1= \prt^-_1 \kappa_2, \; q \kappa_3=
\prt^+_1\kappa_2,\\
q\lambda_1= \prt^-_1 \lambda_2, \; q \lambda_3= \prt^+_1\lambda_1.
\end{gather}

We now use the fact that all homotopies are rel vertices and that
the maps $\alpha, \alpha', \beta, \beta':I^2 \to B$ are constant
on the edges $\prt^-_2,\prt^+_2$. So in the following picture, the
dotted lines represent constant paths, and $\boldsymbol{\zeta}$ is
a hollow cube not yet filled in, but has four faces well defined.
Note also that the maps $I^2 \to B$ given by
$\prt^\pm_2(\kappa_2)$ and $\prt^\pm_2(\lambda_2)$ are constant
maps, by our definition of homotopies.

$$\def\labelstyle{\textstyle}\xymatrixrowsep{1pc}
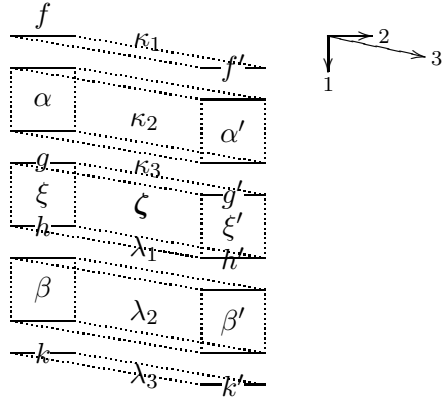
\xymatrix@!0 @M=0pt{{} \ar @{-} [r] ^f  \ar @{.} [drrr]^(0.7){\kappa_1}&\ar @{.} [drrr]&& &\\
\ar @{-} [r] \ar @{.} [dd]\ar @{.} [drrr] \ar @{} [ddr] |\alpha
 &\ar @{.} [drrr] \ar @{.} [dd]
&&\ar @{-} [r]|{f'}&
\\ \ar @{} [drrr] |(0.7){\kappa_2}&&&\ar @{-} [r] \ar @{.} [dd] \ar @{} [ddr] |{\alpha'}  &\ar @{.} [dd]\\
\ar @{-} [r]\ar @{.} [drrr] &\ar @{.} [drrr]&&&\\ \ar @{-} [r] |g
\ar @{} [ddr] |\xi\ar @{.} [drrr]^(0.7){\kappa_3}\ar @{} [drrr]
_(0.7){\boldsymbol{\zeta}}\ar @{.} [dd]&\ar @{.} [dd]\ar @{.}
[drrr]&&\ar @{-} [r]&
 \\
&&& \ar @{-} [r]|{g'}\ar @{.} [dd]\ar @{} [ddr] |{\xi'}&\ar @{.} [dd]\\
\ar @{-} [r]|h \ar @{.} [drrr] \ar @{.} [drrr] |(0.7){\lambda_1}
&\ar @{.} [drrr]  &&&
\\ \ar @{-} [r] \ar @{.} [drrr]\ar @{-} [r]\ar @{.} [dd]\ar @{} [ddr] |\beta &
\ar @{.} [drrr]\ar @{.} [dd]&&\ar @{-} [r]|{h'}\ar @{-} [r]&\\
\ar @{} [drrr] |(0.7){\lambda_2}&&&\ar @{.} [dd] \ar @{-} [r]\ar
@{} [ddr] |{\beta'}&\ar @{.} [dd]\\\ar @{.} [drrr] \ar @{-}
[r]&\ar @{.} [drrr]&&&\\\ar @{.} [drrr] \ar @{-} [r] |k \ar @{}
[drrr] |(0.7){\lambda_3}&\ar @{.} [drrr]&&\ar @{-} [r]&\\ &&&\ar
@{-} [r] |{k'}&\\&&&&}\qquad \threeaxes{2}{1}{3}
$$

The maps $\xi,\xi',\kappa_3,\lambda_1 $ define a map $$ (I^2
\times \dot{I}) \cup (\dot{I} \times I^2) \to E $$ (where
$\dot{I}=\{0,1\}$) given by $ (s,t,0) \mapsto \xi(s,t), (s,t,1)
\mapsto \xi'(s,t)$ on $I^2 \times \dot{I}$, and by $(0,t,u)
\mapsto \kappa_3(t,u), (1,t,u) \mapsto\lambda_1(t,u)$ on
$\dot{I}\times I^2$ respectively. By the rel vertices condition,
these maps can be extended by the constant map over $I \times
\{0\} \times I$, which is $\prt^-_2 (\boldsymbol{\zeta})$. So we
now have maps defined on 5 faces of $I^3$ and agreeing on their
common edges, and so these extend to a map $ \zeta: I^3 \to E$.
However, while this map does agree with the other homotopies, the
result will not be a homotopy of the type required since $\zeta_1=
\prt^+_2(\zeta)=(\zeta|(I \times \{1\}\times I)$ is not constant
as would be required. So we have to make a modification to get a
homotopy between representatives of the original classes.
Intuitively, we move the face $\zeta_1$ of $\zeta$ to the right
and down of our composite picture. This modification will also
change   $\lambda_2, \lambda_3$, but this does not matter for our
purposes, since we need to show only that a homotopy of the
required type exists.

We let $\eps,\eps_1,\eps_2$ denote degenerate elements -- the
element they act on in the following formulae will be clear from
the context, in order to make the compositions properly defined.

Our new homotopy $$ (f\circ \eps,\A\circ_2 \eps_2,g\circ
\eps),(g\circ \eps, ((q\xi)\circ_1 \beta)\circ_2 \eps_2, k\circ
\eps)) \equiv (f'\circ \eps,\A'\circ_2 \eps_2,g'\circ
\eps),(g'\circ \eps, ((q\xi')\circ_1 \beta')\circ_2 \eps_2,
k'\circ \eps)) $$ will be given by
\begin{align}
\kappa_1\circ_2 \eps : f \circ \eps &\simeq  f' \circ \eps, \\
\kappa_2\circ_2 \eps_2 : \alpha \circ_2 \eps_2 & \equiv  \A' \circ_2
\eps_2,\\
\kappa_3\circ_2 \eps : g \circ \eps &\simeq  g' \circ \eps,\\
\lambda'_2: ((p\xi)\circ_1 \beta) \circ_2 \eps_2  & \simeq ((p
\xi') \circ_1 \beta') \circ_2 \eps_2 ,\\
\lambda_3 \circ \prt^+_2 \zeta: k \circ \eps & \simeq k' \circ
\eps ,
\end{align}
where
$$ \lambda_2' = \begin{bmatrix} q \zeta& \tr \;\\
\lambda_2 & \vv\end{bmatrix} \qquad \directs{2}{1}  $$ where
$\tr\;$ is given by $(s,t,u) \mapsto (\prt^+_2\zeta
)(\min(s,1-t),u)$. Note that the combination of $\tr\;$ and
$\eps=\vv\;$  in the second column of the matrix has the effect of
pushing the non constant face $ \prt^+_2\zeta$ of $\zeta$ down to
be able to combine with $\lambda_3$.

It is clear that compositions with degenerate elements in
direction 2 do not change homotopy classes, and so this completes
our geometric proof that $\Phi$ is well defined.

Next we must prove $\Phi\Xi=1, \Xi \Phi =1$.

Considering the formula \eqref{P} for $\Phi$, we see that we can
set $$ \Phi([f,\A,g],[g,\beta,k])= [(f,\A,g),(g,(q\xi)\circ_1
\beta,k)]$$ where now $\xi$ can be chosen to be a constant
homotopy $\eps$. It is then easily seen that $$[(f,\A,g),(g,(q\eps
)\circ_1 \beta,k)]=[(f,\A,g),(g,\beta,k)],$$ and so that
$\Phi\Xi=1$.

To prove $\Xi\Phi=1$ it is sufficient to  show that if $\xi: g
\simeq h$ is a homotopy rel vertices, then $ (g, (q\xi)\circ_1
\beta,k)\equiv( h,\eps_1(q\xi) \circ_1 \beta,k) $. Such a homotopy
is
 given by $(\xi,\kappa, \eps_3(k))$ where $$\kappa= \begin{bmatrix}
 \br \;\\ \eps_3(\beta) \end{bmatrix} \qquad \directs{3}{1}.$$
 \end{proof}
\begin{cor}
The composition $\circ_1$ on $R_2(q)$ is inherited by $\rho_2(q)$
so that $\rho( q)$ becomes a double groupoid.
\end{cor}
\begin{proof}
The composition  $\circ_1$ on $\rho_2( q)$  is the composition of
the maps
$$ \rho_2(q) {\,}_{\prt^+_1}\!\! \times _{\prt^-_1}\rho_2(q)
\labto{\Phi}
             \pi_0^v( R_2(q) {\,}_{\prt^+_1}\!\! \times _{\prt^-_1}R_2(q))
  \to   \rho_2(q)
$$
where the second map is induced by the composition on $R_2(q)$. It
is easy to see that the structure $\circ_1$ gives a groupoid
structure on $\rho_2(q)$.

Thus the only part remaining is the interchange law. However we
easily find that a double composition can be given as
$$ \begin{bmatrix}[f,\A,g] & [f', \A',g'] \\ {}
[h,\beta,k] & [h',\beta',k']
\end{bmatrix} = \left[ f \circ f', \begin{bmatrix}\A & \A'
\\{} (q \xi)\circ_1 \beta & (q \xi')\circ_1
\beta'\end{bmatrix}, k \circ k'\right] \qquad \directs{2}{1}
$$ where $\xi: g \simeq h, \xi': g' \simeq h'$. So the interchange
law follows from that for singular squares.
\end{proof}

Finally, we have to show the relation between the composition
$\circ'$ on $K_1M'$ and the composition $\circ_1$ above. Let
$(f,\A,g),(g,\beta,h) \in R_2(q)$. We first note that, analogously
to the existence of identities in the fundamental groupoid,
$$[f,\A,g] = [f,\A\circ_1 (\eps_1(qg)),g],\; [g,\beta,h]=[g,(\eps_1(qg))\circ_1
\beta,h].$$ Hence
\begin{align*}[f,\A\circ_1\beta,h]=&
[f,\A,g]\circ_1[g,\beta,h]\\=& [f,\A
\circ_1\eps_1(qg),g]\circ_1[g,\eps_1(qg)\circ_1 \beta,h]\\=&
[(f,\A\circ_1 \eps_1(qg),g)\circ'(g,\eps_1(qg)\circ_1
\beta,h)]\end{align*} as required.

\section{Examples}
In order to study the double groupoid $\rho(q)$ we need to have
examples of double groupoids with which to compare it, in addition
to the product fibration of Example \ref{prodfib}. As we shall
see, there are some sub-double groupoids of $\rho(q)$ which are
familiar, but it is interesting that we have little information
about the most general form of double groupoids. For example, the
methods of \cite{BM} give an equivalence between double groupoids
satisfying some filler conditions and what are there called {\em
core diagrams}, but these do not seem to be helpful in this case.

Here we suggest various examples and comparisons for further
investigation.
\begin{example}
Let $i:M \to B$ be the inclusion of a subspace $M$ of $B$. Then
the equivalence relation $\Eq(i)$ is discrete, and so $\rho(i)$ is
a 2-groupoid. Further, if $m \in M$ then the natural map $$ \eta:
\pi_2(B,M,m) \to \rho_2(i)$$ is injective.

\begin{proof}
We represent $\pi_2(B,M,m)$ by maps $\alpha: I^2 \to B$ such that
the face $\pt^-_1 \alpha$ maps into $M$ and the other three faces
map to the base point $m$. The homotopy classes of $\alpha$ which
yield an element $[\alpha]$ of  $\pi_2(B,M,m)$ are through maps of
the same type. Then $\alpha$ also yields an element $\lan \alpha
\ran$ of $\rho_2(i)$, but there the homotopies allow $\pt^+_1
\alpha$ to vary in $M$. We have to prove that the map $\eta:
[\alpha] \mapsto \lan \alpha \ran$ is injective. Suppose then
$[\alpha_-], [\alpha_+] \in \pi_2(B,M,m)$ and $\lan\alpha_-\ran =
\lan\alpha_+\ran \in\rho_2(i)$. Let $h: I^3 \to B$ be a homotopy
determined by this equality, so that
$$\pt^-_3(h)=\alpha_-,\pt^+_3(h)=\alpha_+,$$ $\pt^\pm_2(h)$ maps to
$m$. Let $\theta=\pt^+_1(h)$. The problem is that $\theta$ is not
constant. So we change $h$ to  `move' $\theta$ to the top face and
still give a homotopy $h':\alpha_-'\simeq\alpha_+' $ where
$[\alpha_\pm]=[\alpha_\pm ']$. We can take
\begin{align*} h'&=\begin{bmatrix} \vv & h \\ \bl & \br
\end{bmatrix} \qquad \directs{2}{1}\\ \intertext{so that the two ends of this homotopy are }
\pt^\pm _3(h') &= \begin{bmatrix}\sq & \alpha_\pm\\\sq & \sq
\end{bmatrix},\end{align*} where $\sq$ denotes a double identity, as required.
\end{proof}

Note that $\rho(i)$ is the homotopy 2-groupoid of a pair discussed
by Moerdijk and Svensson in \cite{MandS:2-types}, and is also
recovered from the work of \cite{KP}.

\end{example}

\begin{example}
 The double groupoid $\rho(q)$ contains a 2-groupoid $$ \xymatrix
{ \bar{\rho}(q) \ar@<1ex> [r] \ar@<-1ex> [r] & \pi_1(M) \ar@<1ex>
[r] \ar@<-1ex> [r] & M} $$ where $ \bar{\rho}_2(q)$ is the subset
of   $ \rho_2(q)$ of elements $u$ such that $\prt^-_2u, \;
\prt^+_2u$ are degenerate, that is consist of pairs $(x,x)$. This
is essentially the homotopy 2-groupoid of a map discussed by Kamps
and Porter in \cite{KP}. This 2-groupoid contains various
cat$^1$-groups of the form considered by Loday in \cite{Lo}. The
crossed module  of groupoids associated to this 2-groupoid is of
the form $ C \to \pi_1(M)$ where for each point $x \in M$ we have
$C(x)$ is isomorphic to $\pi_1(F_x, \bar{x})$, the fundamental
group of the homotopy fibre $F_x$ of $q$ over $q(x)$ at the base
point $\bar{x}$ determined by $x$. If $M$ is a subspace of $B$ and
$q$ is the inclusion then $C(x)$ is isomorphic to the familiar
relative homotopy group $\pi_2(B,M,x)$ and the crossed module
$C(x) \to \pi_1(M,x)$ is essentially that first studied by J.H.C.
Whitehead. However we do not  have a reconstruction method for
$\rho(q)$ from  $ \bar{\rho}(q)$, whereas the 2-groupoid can be
reconstructed from the crossed module of groupoids it contains, as
shown in \cite{BH4}.\hfill $\Box$
\end{example}
\begin{example}{\bf Foliations} Let $\cF$ be a foliation on a space
$M$. Thus the leaves of the foliation define an equivalence
relation $R= R(\cF)$. Let $q:M \to B$ be a map of spaces. The
foliation defines a finer topology than that given on $M$ to give
a space $M_\cF$ in which all leaves of the foliation are open
components. So we also have a map $q_\cF : M_\cF \to B$ and hence
may define the homotopy double groupoid $\rho(q_\cF)$. Where this
differs from $\rho(q)$ is that in $\rho(q_\cF)$ the `horizontal'
paths, and the homotopies of paths, all lie in leaves of the
foliation.

An illustrative example is the M\"{o}bius Band $M$ with its
projection $q : M \to S^1$ and foliation $\cF$ by circles of which
the centre one goes once round the Band and the other circles go
twice round. Then $\rho( q_\cF)$ contains the  double groupoid
$\cD(M)$ explained in the Introduction, and which seems to be a
good discrete algebraic model of the foliated M\"{o}bius Band.
\hfill $\Box$
\end{example}

\section*{Acknowledgements}
This work was partially supported by the following grants:INTAS
93-436 `Algebraic K-theory, groups and categories', 97-31961
`Algebraic Homotopy, Galois Theory and Descent', `Algebraic
K-theory, Groups and Algebraic Homotopy Theory'; with Bielefeld,
an  ARC Grant 965 `Global actions and algebraic homotopy', and by
the London Mathematical Society fSU Scheme.

The first author is also grateful to the Erwin Schr\"{o}dinger
Institute of Mathematical Physics and a Leverhulme Emeritus
Fellowship for support to attend a Workshop on Foliations in
August, 2002. 
\end{document}